\documentclass[floatfix,nofootinbib,longbibliography,notitlepage,superscriptaddress]{revtex4-1}

\usepackage{amsmath,amsthm,amssymb,amsfonts,physics,graphicx,mathrsfs}
\usepackage{xcolor}
\usepackage{hyperref}
\usepackage[section]{algorithm}
\usepackage{algpseudocode}
\usepackage{enumitem} 

\newtheorem {theorem} {Theorem}
\newtheorem {result} {Result}

\newtheorem {proposition}{Proposition}
\newtheorem {remark} {Remark}

\newcommand*{\mcal}{\mathcal}
\newcommand*{\mbb}{\mathbb}
\newcommand*{\mrm}{\mathrm}

\newcommand*{\mbf}{\mathbf}

\def\hilb{\mathbb{C}^n}

\def\ii{\mathbb{I}}

\def\prob{\mrm{Prob}}
\def\t{\tilde}
\def\pthresh{p_{\textup{thr.}}}
\def\conf{\delta_{\textup{conf.}}}

\def\phaar{P_{\textup{Haar}}}

\def\floating{\epsilon_0}

\def\haar{\textup{Haar}}

\def\wt{\mrm{wt}}
\def\darg{(\,\cdot\,)}

\def\rr{\mbb{R}}
\def\cc{\mbb{C}}
\def\ee{\mbb{E}}
\def\ff{\mbb{F}}
\def\nn{\mbb{N}}

\newcommand{\anynorm}[2]{\left\lVert#1\right\rVert_{#2}}
\newcommand{\lienorm}[2]{\left\lVert#1\right\rVert_{\lie{#2}}}
\newcommand{\lie}[1]{\mathfrak{#1}}

\begin{document}

\title{Certifying Numerical Decompositions of Compact Group Representations}
\author{Felipe Montealegre-Mora}
\affiliation{Institute of Theoretical Physics, University of Cologne, Germany}
\author{Denis Rosset}
\affiliation{Perimeter Institute of Theoretical Physics, Waterloo, Canada.}
\author{Jean-Daniel Bancal}
\affiliation{Université Paris-Saclay, CEA, CNRS, Institut de Physique Théorique, 91191, Gif-sur-Yvette, France}
\author{David Gross}
\affiliation{Institute of Theoretical Physics, University of Cologne, Germany}

\begin{abstract}
  We present a performant and rigorous algorithm for certifying that a matrix is close to being a projection onto an irreducible subspace of a given group representation.
  This addresses a problem arising when one seeks solutions to semi-definite programs (SDPs) with a group symmetry.
  Indeed, in this context, the dimension of the SDP can be significantly reduced if the irreducible representations of the group action are explicitly known.
  Rigorous numerical algorithms for decomposing a given group representation into irreps are known, but fairly expensive.
  To avoid this performance problem, existing software packages -- e.g.\ RepLAB, which motivated the present work -- use randomized heuristics. 
  While these seem to work well in practice, the problem of to which extent the results can be trusted arises.
  Here, we provide rigorous guarantees applicable to finite and compact groups, as well as a software implementation that can interface with RepLAB.
  Under natural assumptions, a commonly used previous method due to Babai and Friedl runs in time $O(n^5)$ for $n$-dimensional representations.
  In our approach, the complexity of running both the heuristic decomposition and the certification step is $O(\max\{n^3\log n,D\, d^2\log d\})$, where $d$ is the maximum dimension of an irreducible subrepresentation, and $D$ is the time required to multiply elements of the group.
  A reference implementation interfacing with RepLAB is provided.
\end{abstract}

\maketitle

\section{Introduction}
\label{sec:intro}

Semi-definite programming is a widely used numerical tool in science and engineering.
Unfortunately, runtime and memory use of SDP solvers scale poorly with the dimension of the problem. 
To alleviate this issue, symmetries can often be exploited to significantly reduce the dimension~\cite{sym-sdps,jordan-alg,markus,sos,poly-opt,BellIneqSymm,collins2002, BellSymWreath} (see~\cite{invariant-sdps} for a review). 
This requires finding a common block-diagonalization of the matrices representing the symmetry group action.
A large number of numerical methods for this task have been developed \cite{eberly-thesis,eberly-decompositions,mkkk,mm,deklerk,abed-block,cai-algebraic,cai-matrixpoly,cai-livec,replab1,
error-controlled,cai-identification,cai-perturbation,babaiAlgebra}.
These algorithms must be compared along a number of different dimensions:
\begin{enumerate}
  \item
	What is their runtime as a function of the relevant parameters? The most important parameters are the dimension $n$ of the input matrices, the dimension of the algebra $\mathcal{A}$ they span, and the dimension $d$ of the largest irreducible component?
  \item
	Are they probabilistic or deterministic?
  \item
	Do they assume a group structure, or do they work for algebras more generally?
  \item
	Can they handle a situation where only noisy versions of the matrices representating the symmetry are available?
  \item
	Which aspects are covered by rigorous performance guarantees?
\end{enumerate}
While a detailed review of the extensive literature is beyond the scope of this paper, we summarize the performance of the approaches that come closest to the methods described here.

References~\cite{error-controlled,cai-identification,cai-perturbation,babaiAlgebra} give algorithms for finding a block decomposition for general $*$-algebras and come with rigorous guarantees.
Refs.~\cite{cai-identification,cai-perturbation} require one to solve a polynomial optimization problem of degree 4 on $\cc^{n\times n}$.
While this might work in practice, there is no general polynomial-time algorithm for this class of problems.
The procedure of~\cite{error-controlled} requires one to diagonalize ``super-operators'', i.e.\ linear maps acting on $n\times n$-matrices.  
This implies a runtime of $O(n^6)$.

The method of~\cite{babaiAlgebra} exhibits a runtime of $O(\max\{n^2\dim^2\mcal{A},n^3\dim\mcal{A}\})$.
In this scaling, the first term comes from finding an orthogonal basis for $\mcal{A}$ and the second term arises from using this basis to project onto the commutant and to diagonalize.\footnote{This scaling refers to Alg.~B from that reference.
There, the scaling of the second term is presented as $O(n^4\dim\mcal{A})$. 
Upon a closer inspection of their algorithm we found that its runtime is slightly better than claimed.
It seems that the origin of the difference, in their language, is that Alg.~B -- as opposed to Alg.~A -- does not require to use the subroutine \emph{Split}.
Instead, Alg.~B projects a single random matrix onto the commutant of $\mcal{A}$, using $O(n^3\dim\mcal{A})$ operations.
}
While the method comes with a guarantee that the output decomposition is close to invariant, it does not guarantee that the components will be irreducible in the presence of noise.
The runtime is particularily competitive for ``small'' algebras:
If $\alpha\in[0,2]$ is such that $\dim\mathcal{A}=O(n^\alpha)$, the scaling becomes $O(n^{3+\alpha})$ for the case $\alpha < 1$.
On the other hand, in the regime $\alpha>1$, the runtime $O(n^{2+2\alpha})$ is worse than other methods discussed below.

Reference~\cite{babai91} works on finite group representations, rather than general $*$-algebras.
It generalizes Dixon's method~\cite{dixon} to handle noise in the group representation.
This algorithm produces a guaranteed full decomposition, however, for this it must project a full matrix basis onto the commutant of the representation and diagonalize each projection.
This means that its runtime scales quite steeply, as $O(n^5)$, with $O(n^2)$ matrix diagonalizations.

Here, we suggest to split the problem of decomposing a unitary group representation $\rho$ on $\hilb$ into three steps:
\begin{enumerate}
  \item 
  Use a fast heuristic to obtain a candidate decomposition $\hilb\simeq R_1\oplus R_2 \oplus\dots$.
  One particular randomized algorithm running in time $O(n^3)$ has been analyzed 
  \cite{replab1,tavakoli} 
  and implemented as part of the RepLAB \cite{replab-git} software package by some of the present authors.
  While this algorithm seems to give accurate results in practice, this is not underpinned by a formal guarantee.

  \item
  Certify that each of the candidate spaces $R_i$ is within a pre-determined distance $\epsilon$ of a subspace $K_i$ that is invariant under the group.

  \item
  Certify that the invariant spaces $K_i$ are irreducible.
\end{enumerate}

With the first step already covered in Ref.~\cite{replab1,tavakoli}, the present paper focuses on the two certification steps.
Thus, we are faced with the situation that a heuristically obtained $n\times n$ matrix $\pi$ is provided, which may or may not be close to a projection onto an invariant and irreducible space.
We provide a probabilistc algorithm for this decision problem.
More precisely, our main result is this:

\begin{result}
  Let $G$ be a compact group.
  Assume that:
  \begin{enumerate}
	\item
	  There exists a representation $g\mapsto \rho(g)$ in terms of unitary $n\times n$ matrices.
	\item
	  In time $O(n^2)$, one can draw an element $g\in G$ according to the Haar measure, and compute an approximation $\tilde\rho$ such that $\max_{g}\max_{ij} |\rho_{ij}(g)-\tilde\rho_{ij}(g)| = o\left(\frac1{n^3\log n}\right)$.
  \end{enumerate}
  Then there exists an algorithm that takes as input an $n\times n$ matrix $\pi$ as well as numbers $\epsilon, \pthresh$, and returning $\mathrm{true}$ or $\mathrm{false}$ such that:
  \begin{enumerate}
	\item{}[False positive rate]
	  The probability that the algorithm returns $\mathrm{true}$ even though $\pi$ is not $\epsilon$-close in Frobenius norm to a projection onto an invariant and irreducible $\rho$-space is \emph{upper-bounded} by $\pthresh$.
	\item{}[False negative rate]
	  The probability that the algorithm returns $\mathrm{false}$ even though $\pi$ is $(\epsilon/2)$-close in Frobenius norm to a projection onto an invariant and irreducible $\rho$-space is \emph{approximately} $2\pthresh$.
	\item{}[Runtime]
	  As long as $\epsilon = o\left(\frac1{n^2\log n}\right)$, the algorithm terminates in time
	  \begin{align*}
		O\left(\big(n^3\log n+D\tr(\pi)^2\log\tr\pi\big)\log\frac1{\pthresh}\right),
      \end{align*}
	  where $D$ is time required to multiply two elements of $G$.
  \end{enumerate}
\end{result}

This algorithm has been implemented in Python and is available in~\cite{repcert}.

There is an asymmetry in the way we treat false positives rates (which are bounded rigorously) and false negative rates (which are only approximated). 
This reflects the different roles these two parameters play in practice.
Indeed, if the certification algorithm returns $\mathrm{false}$, the symmetry reduction has failed, no further processing will take place, and thus no further guarantees are needed.
In contrast, if the algorithm returns $\mathrm{true}$, the user must be able to quantify their confidence in the result -- hence the necessity to have a rigorous upper bound on the false positive rate.

In the main text, we introduce a an additional parameter $\delta$, which can be used to tune the false negative rate independently of the false positive rate $\pthresh$.
The interpretation is that $\delta$ is a rigorous upper bound on the false negative rate in the limiting case where $\epsilon=0$ and the approximation $\tilde\rho$ is in fact exact.
We have chosen $\delta = 2\pthresh$ in the displayed result, which turns out to simplify the formula for the runtime.

In practice, one can find appropriate values for $\delta$ numerically: In an \emph{exploratory phase}, one can run the algorithm for increasing values of $\delta$, until it reliably identifies valid inputs as such.
One would then certify a subspace by running the procedure \emph{once} with the $\delta$ previously obtained.

The paper is organized as follows.
In Sec.~\ref{sec:prelim} we review the mathematical setting of the paper.
In Sec.~\ref{sec:invariance} and Sec.~\ref{sec:irreducibility} we present the algorithms to certify invariance and irreducibiltity respectively.
Finally, in Sec.~\ref{sec:time} we discuss the runtime of the algorithms.

\section{Mathematical setting}
\label{sec:prelim}

Let $G$ be a compact group, and $(\hilb,\rho)$ be a unitary representation of $G$. 
A subset $S\subset G$ \emph{generates} the group if $\langle S\rangle$ is dense in $G$, and it is \emph{symmetric} if $S=S^{-1}$.

We assume that the user can evaluate a function $\t\rho: G\to \cc^{n\times n}$ satifying
\begin{align*}
  \max_{ij}|\rho(g)_{ij} - \t\rho(g)_{ij}| \leq \floating, \qquad \forall\,g\in G.
\end{align*}
If $R\subset\hilb$ is the subspace to be certified and $\pi_R$ projects onto it, we use $\t\pi_R$ to denote an approximation to $\pi_R$:
\begin{align*}
  \max_{ij}|(\pi_R)_{ij} - (\t\pi_R)_{ij}| \leq \floating.
\end{align*}
We require that $\floating<\frac 1{2n}$, however in practice $\floating$ is typically of the order of machine precision.

In the context of our algorithms, the user has obtained $\t\pi_R$ as an output of their
numerical procedure to decompose $\rho$. 
Using this operator as an input, the goal is to certify two statments. 
The first is that there exists some invariant subspace $K\subset\hilb$ with associated projector $\pi_K$ satisfying that 
\begin{align}
  \label{eq:approx invariance}
  \|\pi_R - \pi_K\|_F \leq \epsilon,
\end{align}
where $\|\cdot\|_F$ is the Frobenius norm and the precision parameter $\epsilon<1/2$ is an input.
We call this procedure \emph{certifying invariance}.
The second is that the subspace $K$ is an irreducible $G$ representation. 

For this task, we assume that one \emph{1.}\ knows an upper bound $r_G$ on the number of generators of $G$, and \emph{2.}~can sample from the Haar measure and evaluate $\t\rho$  on the sample.
In an appendix, we show how to relax the second condition and instead assume only that the user can evaluate $\t\rho$ on a well-behaved \emph{fixed} generator set.
The algorithms are probabilistic.
A bound $\pthresh$ on the false positive rate -- i.e.\ the probability that an input is certified even though it is not close to the projection onto an irredudcible representation -- is an explicit parameter.

Bounds $r_G$ on the number of generators of $G$ are known for a wide variety of groups.
For example it is known that $r_G\leq2$ when $G$ is a finite dimensional connected compact group~\cite{connectedRank}. 
For a wide variety of finite simple groups, furthermore, $r_G\leq7$ (see~\cite{simpleRank} for a review).


\section{The invariance certificate}
\label{sec:invariance}

Here we present our algorithm for the first task, that is, certifying the approximate invariance of $R$.
Section~\ref{sec:closeness} treats a closely related problem: deciding whether an operator is close to the \emph{commutant}
\begin{align*}
  \{ Y \in\mathbb{C}^{n\times n}\,|\, [\rho(g), Y ] = 0 \;\forall g\,\in G \}
\end{align*}
of $\rho$.
In that section we also work in the idealized case where $\floating=0$.
The general algorithm deciding invariance is presented in Section~\ref{sec:inv}.

\subsection{Estimating closeness to the commutant in the ideal case}
\label{sec:closeness}

As mentioned, in this section we assume $\floating=0$ -- i.e.\ that the representation $\rho$ can be evaluated \emph{exactly} -- in order to bring out the key components of the argument.

Consider an $n\times n$ matrix $X$ (later, we will take $X$ to be the approximate projection $\t\pi_R$ onto a candidate subspace).
The randomized Algorithm~\ref{alg:closeness} tests whether 
\begin{align*}
  \|X - \phaar(X)\|_\infty \leq \epsilon.
\end{align*}
There, $\|\cdot\|_\infty$ is the spectral norm and $\phaar$ is the Hilbert-Schmidt projection onto the commutant
\begin{align*}
  \phaar(X) := \ee_g[\rho(g)X\rho^\dagger(g)],
\end{align*}
where the expectation value is with respect to the Haar distribution. 

\begin{algorithm}[H]
  \caption{Closeness to Commutant}
  \label{alg:closeness}
  \hspace*{\algorithmicindent} \textbf{Input:}%
  \begin{itemize}
    \item $X \in \cc^{n\times n}$,
	\item $\pthresh \in (0,1), \; \epsilon \in (0,1/2)$.
  \end{itemize}
  \begin{algorithmic}[1]
    \State Set $r=8\lceil(\log(1/\pthresh) + \log(2n))\rceil$\;
    \vspace{1mm}
    \State Sample $r$ group elements $g_1,\dots,g_r\in G$ Haar-randomly\;
    \vspace{1mm}
    \State Compute $c=\Big\|\frac1r \sum_i\rho(g_i)X\rho^\dagger(g_i) - X\Big\|_\infty$\;
    \If{$2c\leq\epsilon$}
      \State \textbf{Return:} \emph{True}\;
    \EndIf
    \State \textbf{Return:} \emph{False}\;
    %
  \end{algorithmic}
\end{algorithm}

\begin{proposition}
  \label{prop:closeness}
  Let $X\in\cc^{n\times n}$ satisfy $\|X-\phaar(X)\|_\infty > \epsilon$. 
  Then, the probability that Alg.~\ref{alg:closeness} returns \emph{True} is at most $\pthresh$.
\end{proposition}
\begin{proof}
  Consider the following matrix-valued random variable with mean equal to zero,
  \begin{align*}
    Z_g := \frac1r\Big( \rho(g)X\rho^\dagger(g) - \phaar(X) \Big),
    \qquad
    g\in G \text{ Haar random}.
  \end{align*}
  Using $R := \mrm{Id} - \phaar$ (the projector onto the orthocomplement of the commutant of $\rho$), we find $Z_g=\frac1r \rho(g) R(X) \rho^\dagger(g)$, and so,
  \begin{align*}
    \|Z_gZ_g^\dagger\|_\infty
    =
    \frac{1}{r^2}\|R(X)R(X)^\dagger\|_\infty
    =
    \frac{1}{r^2}\|R(X)\|_\infty^2, \qquad \forall\,g\in G. 
  \end{align*}
  
  This way, by the matrix Hoeffding bound~\cite{mackey},
  \begin{align*}
    \prob\left[
      \Big\| \sum_i Z_{g_i} \Big\|_\infty \geq z \| R(X)\|_\infty
    \right]
    \leq
    2n\exp(\frac{-rz^2}{2})
  \end{align*}
  where $\{g_i\}$ are the samples in line~2 of Alg.~\ref{alg:closeness}.
  Taking $z=1/2$, the right-hand side above is $\leq\pthresh$ and so with probability at least $1-\pthresh$ it holds that
  \begin{align*}
    c = 
    \Big\| \frac1r \sum_i \rho(g_i)X\rho^\dagger(g_i) - X \Big\|_\infty
    =
    \Big\| \sum_i Z_{g_i} - R(X) \Big\|_\infty
    \geq
    \| R(X) \| - \Big\| \sum_i Z_{g_i} \Big\|_\infty
    \geq
    \frac12 \| R(X)\|_\infty > \epsilon/2.
  \end{align*}
\end{proof}

We now show a converse result, namely, that Alg.~\ref{alg:closeness} always “detects" matrices which are close enough to the commutant.

\begin{proposition}
  \label{prop:closeness converse}
  Let $X$ satisfy $\|X-\phaar(X)\|_\infty\leq \epsilon/2$ for some $\epsilon<1$.
  Then Alg.~\ref{alg:closeness} deterministically returns \emph{True} upon the input $X$, $\epsilon$.
\end{proposition}
\begin{proof}
  For any $g\in G$ it holds that
  \begin{align*}
    \|[\rho(g),X]\|_\infty = 
    \|[\rho(g),X-\phaar(X)]\|_\infty \leq 
    2\|X-\phaar(X)\|_\infty \leq 
    \epsilon.
  \end{align*}
  Therefore, using standard norm relations we obtain 
  \begin{align*}
    c = 
    \Big\| \frac1r \sum_i \big(\rho(g_i)X\rho^\dagger(g_i) - X \big) \Big\|_\infty
    \leq
    \frac1r \sum_i 
    \Big\| [\rho(g_i),X] \Big\|_\infty \leq \epsilon.
  \end{align*}

\end{proof}

\subsection{The full certificate}
\label{sec:inv}

Here, we will go beyond Section~\ref{sec:closeness} in two ways:
First, we allow for non-zero errors $\epsilon_0$.
Second, we show that a projection that is close to being invariant is close to a projection onto an invariant subspace.
The goal is, given $\t\pi_R$ as an input, to certify that there is an invariant subspace $K$ with
\begin{align*}
  \|\pi_K-\pi_R\|_F \leq \epsilon.
\end{align*}
The procedure is given in Alg.~\ref{alg:invariance}.

\begin{algorithm}[H]
  \caption{Invariance certificate}
  \label{alg:invariance}
  \hspace*{\algorithmicindent} \textbf{Input:}%
  \begin{itemize}
    \item $\t\pi_R\in\cc^{n\times n}$,
    \item $\pthresh\in(0,1)$,
    \item $\epsilon \in (0,1/2)$.
  \end{itemize}
  \hspace*{\algorithmicindent} \textbf{Output:}%
  \emph{True/False}\\
  \hspace*{\algorithmicindent}
  \begin{algorithmic}[1]
    \State Set $r=8\lceil(\log(1/\pthresh) + \log(2n))\rceil$, $f_{\text{err}} = 8n\floating + 6n^2\floating^2 + 2n^3\floating^3$, and $\epsilon'=\epsilon/2\sqrt{2\dim R}$\;
    \vspace{1mm}
    \State Sample $r$ group elements $g_1,\dots,g_r\in G$ Haar-randomly\;
    \vspace{1mm}
    \State Compute $\t c=\Big\|\frac1r \sum_i\t\rho(g_i)\t\pi_R\t\rho^\dagger(g_i) - \t\pi_R\Big\|_\infty$\;
    \vspace{1mm}
    \If{$2\t{c}+f_\text{err}\leq\epsilon'$}
      \State \textbf{Return:} \emph{True}\;
    \EndIf
    \State \textbf{Return:} \emph{False}\;
    %
  \end{algorithmic}
\end{algorithm}

As before, line 4 of Alg.~\ref{alg:invariance} simply takes $k$ close to the minimum of $f_k(c)$ and does not affect the probability of falsely certifying $R$.
Our main result in this section is the following guarantee on the invariance certificate.

\begin{theorem}
  \label{thm:invariance}
  Assume that for all invariant subspaces $K\subset\cc^n$,
  \begin{align}
    \label{eq:not invariant}
    \|\pi_K - \pi_R\|_F > \epsilon.
  \end{align}
  Then, the probability that Alg.~\ref{alg:invariance} returns \emph{True} is upper bounded by $\pthresh$.
\end{theorem}

To prove Thm.~\ref{thm:invariance} we will first show that if $\pi_R$ is close to the commutant, then it is close to an invariant projector $\pi_K$ as in eq.~\eqref{eq:approx invariance}.
After that, our argument will closely follow Sec.~\ref{sec:closeness}.

\begin{proposition}
  \label{prop:far from commutant}
  Assume that $\pi_R$ satisfies $2\sqrt{2\dim R}\,\|\phaar(\pi_R)-\pi_R\|_\infty \leq\epsilon$ for some $\epsilon<1$. 
  Then there exists an invariant subspace $K$ with projector $\pi_K$ satisfying $\|\pi_R-\pi_K\|_F \leq \epsilon$.
\end{proposition}
\begin{proof}
  Let $\lambda^\downarrow(M)$ be the vector of eigenvalues of a Hermitian matrix $M\in\cc^{n\times n}$ in decreasing order.
  By Weyl's perturbation theorem (see e.g.~\cite[Chap. VI]{bhatia}),
  \begin{align*} 
    \|\lambda^\downarrow(\phaar(\pi_R))-\lambda^\downarrow(\pi_R)\|_{\ell_\infty}
    \leq
    \frac{\epsilon}{2\sqrt{2\dim R}} = \epsilon'.
  \end{align*}
  This way, the eigenvalues of $\phaar(\pi_R)$ lie in $[-\epsilon',\epsilon']\cup[1-\epsilon',1+\epsilon']$, where $\epsilon'<1/2$.
  Let $\pi_K$ be the projector onto all eigenspaces corresponding to eigenvalues in $1\pm\epsilon'$. 
  The projector $\pi_K$ is invariant and satisfies $\|\pi_K-\phaar(\pi_R)\|_\infty\leq\epsilon'$. 
  We therefore see that, 
  \begin{align*}
    \|\pi_K - \pi_R\|_F
    &\leq
    \sqrt{2\dim R}\,\|\pi_K - \pi_R\|_\infty\\
    &\leq
    \sqrt{2\dim R}\,\big(
      \|\pi_K - \phaar(\pi_R)\|_\infty+\|\phaar(\pi_R) - \pi_R\|_\infty
    \big)\\
    &\leq
    2\epsilon'\sqrt{2\dim R} = \epsilon,
  \end{align*}
  where we used that $\rank(\pi_K-\pi_R)\leq\dim K + \dim R = 2\dim R$ in the first step.
\end{proof}

From the proof above it becomes clear that certifying that $R$ is approximately invariant is, ultimately, just certifying that $\pi_R$ is close enough to the commutant.

\begin{proof}[Proof of Thm.~\ref{thm:invariance}]
  By Prop.~\ref{prop:far from commutant} we may take
  \begin{align*}
    \frac{\epsilon}{2\sqrt{2\dim R}} < \,\|P_\haar(\pi_R) - \pi_R\|_\infty.
  \end{align*}
  Let
  \begin{align*}
    A :=
    \frac1r\sum_i\Big( 
      \rho(g_i) \pi_R \rho^\dagger(g_i) - \t\rho(g_i) \t\pi_R \t\rho^\dagger(g_i)
    \Big),
    \qquad
    \Delta_R:=\pi_R-\t\pi_R,
  \end{align*}
  then,
  \begin{align*}
    \Big\|
        \frac1r\sum_i\rho(g_i) \pi_R \rho^\dagger(g_i) - \pi_R
    \Big\|_\infty
    \leq
    \|\Delta_R\|_\infty + \| A\|_\infty + 
    \Big\|
      \frac1r\sum_i \t\rho(g_i) \t\pi_R \t\rho^\dagger(g_i) - \t\pi_R
    \Big\|_\infty
    =
    n\floating + \| A\|_\infty + \t{c}.
  \end{align*}
  Then, by Prop.~\ref{prop:closeness}, with probability at least $1-\pthresh$ it holds that 
  \begin{align*}
    \frac{\epsilon}{2\sqrt{2\dim R}} < 2(n\floating + \| A\|_\infty + \t{c}).
  \end{align*}
  We now provide an upper bound on $\| A\|_\infty$.
  Let $\Delta(g):=\rho(g) - \t\rho(g)$, then 
  \begin{align*}
    \|A\|_\infty \leq
    \ee_i\Big[&
      \|\Delta(g_i)\pi_R\rho^\dagger(g_i)\|_\infty + 
      \|\rho(g_i)\Delta_R\rho^\dagger(g_i)\|_\infty +
      \|\rho(g_i)\pi_R\Delta^\dagger(g_i)\|_\infty \\
      + &\|\Delta(g_i)\Delta_R\rho^\dagger(g_i)\|_\infty +
      \|\Delta(g_i)\pi_R\Delta^\dagger(g_i)\|_\infty +
      \|\rho(g_i)\Delta_R\Delta^\dagger(g_i)\|_\infty \\
      + &\|\Delta(g_i)\Delta_R\Delta^\dagger(g_i)\|_\infty
    \Big].
  \end{align*}
  Submultipliciativity, together with $\max\{\|\Delta_R\|_\infty,\,\|\Delta(g)\|_\infty\}\leq n\floating$ for all $g\in G$, gives
  \begin{align*}
    \| A\|_\infty \leq 3(n\floating + n^2\floating^2) + n^3\floating^3.
  \end{align*}

\end{proof}


\section{Irreducibility certificate}
\label{sec:irreducibility}

In this section we present an algorithm that certifies irreducibility.
Given $\t\pi_R$ as an input, where $R$ holds an invariance certificate, the goal is to certify that the minimizer of
\begin{align}
  \label{eq:minimal invariant}
  \min_{\substack{K\subset \hilb\\ K\text{ invar.}}} \|\pi_R - \pi_K\|_F
\end{align}
is irreducible.
We first present the idea of the algorithm in an idealized setting, and then come back to the noisy scenario.

\subsection{The ideal case}
\label{sec:ideal irreducibility}

Let $(\cc^{n_K},\rho_K)$ be a unitary representation of $G$ and suppose that we have access to the same primitives as in Sec~\ref{sec:closeness}. 
Namely, we can sample Haar-randomly from $G$ and evaluate $\rho_K$ on any sample.
Our task is to certify if $\rho_K$ is irreducible. 
The following algorithm uses random walks to acheive this.

\begin{algorithm}[H]
  \caption{Ideal irreducibility certificate}
  \label{alg:ideal irrep}
  \hspace*{\algorithmicindent} \textbf{Input:}%
  \begin{itemize}
    \item $\pthresh\in(0,1)$,\Comment{Bound on false positive rate.}
    \item $\pthresh'\in(\pthresh,1)$,\Comment{Bound on false negative rate.}
  \end{itemize}
  \hspace*{\algorithmicindent} \textbf{Output:}%
  \emph{True/False}.\\
  \hspace*{\algorithmicindent}
  \begin{algorithmic}[1]
    \State Set $r=\max\{r_G,\, 8\lceil(\log(2/\pthresh) + 2\log(n_K))\rceil\}$ \Comment{$G$ generated by $\leq r_G$ elements}\;
    \vspace{1mm}
    \State Set $m=2n_K^2\cdot\max\{8\lceil\log((\pthresh'-\pthresh)^{-1})\rceil,\lceil\log(\pthresh^{-1})\rceil\}$ \Comment{$m$ number of random walks}\;
    \vspace{1mm}
    \State Set $t=2+\lceil\log_2 n_K\rceil$ \Comment{$2t$ length of random walks}\;
    \vspace{1mm}
    \State Sample $r$ elements $g_i\in G$ Haar-randomly and set $S=\{g_i\}\cup\{g_i^{-1}\}$\;
    \vspace{1mm}
    \State Sample $m$ elements $\mbf{s}_i\in S^{2t}$ uniformly\;
    \vspace{1mm}
    \State Compute $E_m=\frac1m\sum_i|\tr\rho_K(\mathbf{s}_i)|^2$\;
    \vspace{1mm}
    \State Set $\theta_m =  n_K\sqrt{2/m}\log(1/\pthresh) $\;
    \vspace{1mm}
    \If{$E_m < 2(1-\theta_m)$}\;
      \State \Return \emph{True}\;
    \EndIf\;
    \State \Return \emph{False}
  \end{algorithmic}
\end{algorithm}

\begin{theorem}
  \label{thm:ideal irrep certificate}
  Let $\rho_K$ be \emph{reducible}, then the probability that Alg.~\ref{alg:ideal irrep} returns \emph{True} upon this input is at most $\pthresh$.
\end{theorem}

Our proof of Thm.~\ref{thm:ideal irrep certificate} will work for any value of $t$, i.e.\ it does not rely on using $t=2+\lceil\log_2n_k\rceil$.
However, if $t$ is chosen too small, the algorithm could fail to recognize irreducible representations ---its false negative rate would be large.
We will bound this rate at the end of this subsection.

The key for the proof of Thm.~\ref{thm:ideal irrep certificate} is Schur's lemma ---if $\rho_K$ were irreducible it would hold that $\tr\phaar=1$ and otherwise it holds that $\tr\phaar\geq 2$.
What the algorithm does is estimate a quantity which is larger than the dimension of the commutant of $\rho_K$. 
As we will see, if $\rho_K$ is reducible then it is exceedingly unlikely for this estimator to fall too much below 2.

The quantity being estimated is, in fact, $\tr P_S^{2t}$, where $P_S$ is the random walk operator associated to $\rho_K$. 
The connection to the dimension of the commutant is made by the following statement.

\begin{proposition}
  \label{prop:ideal upper bound}
  For any $t$ it holds that $\tr\phaar \leq \tr P_S^{2t}$. 
\end{proposition}
\begin{proof}
  Unitarity ensures that $\|P_S\|_\infty=1$. 
  Because $r\geq r_G$, the probability that $S$ generates $G$ is one. 
  Together with $S=S^{-1}$, this ensures that $P_S$ is self-adjoint and that the $+1$ eigenspace corresponds exactly to the commutant of $\rho_K$. 
   
  Let $\{\lambda_i\}$ be all the eigenvalues of $P_S$ different from one. 
  The statement follows from
  \begin{align*}
    \tr P_S^{2t} = \tr\phaar + \sum_i \lambda_i^{2t} \geq \tr\phaar.
  \end{align*}
\end{proof}

\begin{proof}[Proof of Thm.~\ref{thm:ideal irrep certificate}]
  It is clear that $E_m$ is an estimator for $\tr P_S^{2t}$.
  Since $\rho_K$ is unitary, furthermore, $|\tr\rho_K(g)|^2\leq n_K^2$ for any $g$, and so by Chernoff's bound,
  \begin{align*}
    \mrm{Pr}\left[
      E_m \leq (1-\theta)\tr P_S^{2t}
    \right]
    \leq
    \exp(\frac{-\theta^2 m\tr P_S^{2t}}{2n_K^2}),
  \end{align*}
  for any $\theta\in(0,1)$.
  But by the assumption on $m$ we may use $\theta=\theta_m$ in the equation above. 
  Then, using Prop.~\ref{prop:ideal upper bound} and $\tr\phaar\geq2$,
  \begin{align*}
    \mrm{Pr}\left[
      E_m \leq 2(1-\theta_m)
    \right]
    \leq
    \mrm{Pr}\left[
      E_m \leq (1-\theta_m)\tr P_S^{2t}
    \right]
    \leq
    \exp(\frac{-\theta_m^2 m\tr P_S^{2t}}{2n_K^2})
    \leq
    \exp(\frac{-\theta_m^2 m}{n_K^2}) < \pthresh.
  \end{align*}
\end{proof}

As mentioned, the proof above doesn't rely on the particular choice of $t$ in line~3 of Alg.~\ref{alg:ideal irrep}.
It also only uses the bound $m>2n_K^2\log(1/\pthresh)$ on the number of samples (cf. line~2).
In Prop.~\ref{prop:irrep converse}, we use $t>2+\log_2 n$ and $m>16n_K^2\log_2(1/(\pthresh'-\pthresh))$ to bound the false negative rate of the algorithm.
To prove it, it's convenient to show the following intermediate result first.

\begin{proposition}
  \label{prop:bernstein}
  Let $S$ be sampled as in Alg.~\ref{alg:ideal irrep}.
  The probability that $\|\phaar - P_S\|_\infty > 1/2$ is strictly less than
  \begin{align*}
    2n^2\exp\left(\frac{-r}{8}\right)\leq\pthresh.
  \end{align*}
\end{proposition}
\begin{proof}
  Let $\sigma$ be the representation of $G$ acting by conjugation on $\cc^{n\times n}$.
  For a group element $g\in G$ sampled Haar-randomly, the operator
  \begin{align*}
    V_g := \frac1r\Big(\frac12\big(\sigma(g)+\sigma^\dagger(g)\big) - \phaar\Big)
  \end{align*}
  is a Hermitian random variable with zero mean. 
  Furthermore, by unitarity of $\rho$ and because $\sigma(g)$ and $\phaar$ are simultaneously diagonalizable, we have that
  \begin{align*}
    \| V_g \|_\infty \leq \frac1r, \qquad \| V_g^2\|_\infty \leq\frac{1}{r^2}.
  \end{align*}
  But then, writing $S=\{g_i\}_{i=1}^r \cup \{g_i^{-1}\}_{i=1}^r$, we see that
  \begin{align*}
    P_S - \phaar = \sum_{i=1}^r V_{g_i},
  \end{align*}
  where the operators $V_{g_i}$ are independent random variables satisfying the conditions above.
  Then, by the matrix Hoeffding bound~\cite{mackey},
  \begin{align*}
    \prob\left(\lambda_{\max}(P_S - \phaar) > x \right) < n^2 e^{\frac{-rx^2}{2}},
  \end{align*}
  where $\lambda_{\max}$ is the maximum eigenvalue. 
  Finally, repeating the statement above for $\lambda_{\max}(\phaar-P_S)$ and using the union bound, we conclude that
  \begin{align*}
    \prob\left(\|\phaar - P_S\|_\infty > x\right) < 2n^2 e^{\frac{-rx^2}{2}}.
  \end{align*}
  Using $x=1/2$ and the fact that $r\geq8\lceil(\log(1/\pthresh)+ 2\log(n))\rceil$ we recover the claimed statement.
\end{proof}

\begin{proposition}
  \label{prop:irrep converse}
  Let $\rho_K$ be \emph{irreducible}, then the probability that Alg.~\ref{alg:ideal irrep} returns \emph{False} upon this input is at most $\pthresh'$.
\end{proposition}
\begin{proof}
  By Prop.~\ref{prop:bernstein}, with probability at least $1-\pthresh$ it holds that
  \begin{align}
    \label{eq:aux 1}
    \|P_S^{2t}-\phaar\|_\infty \leq 2^{-2t},
  \end{align}
  where we used $P_S^{2t}-\phaar=(P_S-\phaar)^{2t}$ because $P_S$ and $\phaar$ commute.
  This way,
  \begin{align*}
    \tr P_S^{2t} \leq \tr \phaar + n_K^2 2^{-2t} \leq \tr \phaar + \frac{1}{16} = \frac{17}{16}.
  \end{align*}  
  Furthermore, by our assumption in $m$, we have $2(1-\theta_m)\geq3/2$.
  But then, the Chernoff bound says that the probability that $E_m\geq3/2$ is at most
  \begin{align*}
    \exp(\frac{-m}{n_K^2}\frac{49}{3\times256}) < \exp(\frac{-m}{16n_K^2})\leq\pthresh'-\pthresh.
  \end{align*}
  
  A false positive can occur if either eq.~\eqref{eq:aux 1} does not hold, or if conditioned on it holding, $E_m\geq 3/2$.
  By the union bound, this probability is at most $\pthresh+(1-\pthresh)(\pthresh'-\pthresh)<\pthresh'$.
\end{proof}

\subsection{The noisy case}
\label{sec:noisy irrep}

In this section we adapt the idea presented above to the noisy scenario. 
Suppose we have certified that a subspace $R\subset \hilb$ is invariant (with precision $\epsilon$). 
We now wish to certify that the minimizer $K$ of~\eqref{eq:minimal invariant} is irreducible.

The algorithm for this is Alg.~\ref{alg:noisy irrep}.
As before, the algorithm has a controllable false positive rate $\pthresh$ as an input.
This is important from the point of view of certification ---if the output is \emph{True}, then one can be rather certain that $K$ is irreducible.

Additionally, the algorithm takes as an input a confidence parameter $\pthresh<\conf<1$
which roughly tunes the false negative rate.
In fact, this parameter is used in the same way that $\pthresh'$ was used in Alg.~\ref{alg:ideal irrep}.
Because Alg.~\ref{alg:noisy irrep} reduces to Alg.~\ref{alg:ideal irrep} in the limit of  $\epsilon,\,\floating\to0$, we expect that the false negative rate is well approximated by $\conf$ when $\epsilon$ and $\floating$ are small enough.
Since the runtime of the algorithm scales with 
$\max \log(1/\pthresh),\log(1/(\conf-\pthresh))$, 
a reasonable choice for the confidence parameter is $\conf=2\pthresh$.

Within Alg.~\ref{alg:noisy irrep} and throughout this section we use the following conventions:
\begin{align*}
  c_1 :&= 2(\epsilon + n\floating)(1+\epsilon+n\floating) + n\floating(1+\epsilon+n\floating)^2,\\
  c_2 :&= 2c_1(1+c_1),\\
  h_t(x) :&= (1+x)^t - 1,\\
  d_t :&= h_t(c_2),\\
  e_t :&= d_{2t}(\mrm{int}(\tr\t\pi_R)^2 + d_{2t}).
\end{align*}

For the sake of clarity, we have shifted the proofs of several propositions in this subsection to App.~\ref{app:proofs}.

\begin{algorithm}[H]
  \caption{Irreducibility certificate}
  \label{alg:noisy irrep}
  \hspace*{\algorithmicindent} \textbf{Input:}%
  \begin{itemize}
    \item $\t\pi_R\in\cc^{n\times n}$, $\epsilon\in(0,1/2)$ \Comment{$\pi_R$, $\epsilon$ satisfy~\eqref{eq:approx invariance}}
    \item $\pthresh\in(0,1)$ \Comment{Bound on false positive rate}
    \item $\conf$ \Comment{Confidence parameter}
  \end{itemize}
  \hspace*{\algorithmicindent} \textbf{Output:}%
  \emph{True/False}.\\
  \hspace*{\algorithmicindent}
  \begin{algorithmic}[1]
    \If{$e_t\geq 2$}
      \State \Return \emph{False}
    \EndIf 
    \State Set $r=\max\{r_G,\,12\lceil(\log(2/\pthresh) + 2\log(n))\rceil\}$ \Comment{$G$ generated by $\leq r_G$ elements.}\;
    \vspace{1mm}
    \State Set $m = 2\left\lceil\frac{\mrm{int}(\tr\t\pi_R)^2 + d_{2t}}{2-e_t}  \cdot\max\{\log(\pthresh^{-1}),8\log((\conf-\pthresh)^{-1})\}\right\rceil$ \Comment{$m$ random walks}\;
    \vspace{1mm}
    \State Set $t=2+\lceil\log_2 \mrm{int}(\tr\t\pi_R)\rceil$ \Comment{$2t$ random walk length}\;
    \vspace{1mm}
    \State Sample $r$ elements $g_i\in G$, set $S=\{g_i\}\cup\{g_i^{-1}\}$\;
    \vspace{1mm}
    \State Sample $m$ words $\mbf{s}_i\in S^{2t}$ uniformly\;
    \vspace{1mm}
    \State Compute $E = e_t + \frac1m \sum_i|\tr \t\rho_R(\mathbf{s}_i)|^2$\;
    \vspace{1mm}
    \State Set $\theta_m = \sqrt{ 2\log(1/\pthresh) (\mrm{int}(\tr\t\pi_R)^2+d_{2t})/m(2-e_t) }$\;
    \vspace{1mm}
    \If{$E<2(1-\theta_m)$}
      \State \Return \emph{True}
    \EndIf
    \State \Return \emph{False}
  \end{algorithmic}
\end{algorithm}

\begin{theorem}
  \label{thm:irrep certificate}
  Assume that the minimizer $K$ of eq.~\eqref{eq:minimal invariant} is \emph{reducible}. 
  Then the probability that Alg.~\ref{alg:noisy irrep} outputs \emph{True} is at most $\pthresh$.
\end{theorem}

Similar to the ideal case, the proof of this theorem relies on characterizing the \emph{approximate} random walk operator $Q_S^R$ given by
\begin{align*}
  Q_S^R \darg := 
  \frac{1}{|S|}\sum_{s\in S} 
  \t\pi_R\t\rho(s)\t\pi_R^\dagger \darg \t\pi_R\t\rho^\dagger(s)\t\pi_R^\dagger.
\end{align*}
Our approach uses $Q_S^R$ to upper-bound the dimension of the commutant of $\rho$ restricted to $K$, that is $\tr\phaar^K$, where
\begin{align*}
  \phaar^K \darg := \int_G\mrm{d}\mu_\haar(g)
  \pi_K\rho(g)\pi_K\darg \pi_K\rho^\dagger(g)\pi_K.
\end{align*}

An important object in our proof is the \emph{restricted random walk operator},
\begin{align*}
  P_S^K \darg := 
  \frac{1}{|S|}\sum_{s\in S} 
  \pi_K\rho(s)\pi_K\darg \pi_K\rho^\dagger(s)\pi_K.
\end{align*}
Notice that $Q_S^R$ is a small perturbation of $P_S^K$.

\begin{proposition}
  \label{prop:noisy over estimator}
  Use the notation above, let $Q_e := P_S^K - Q_S^R$ and $\gamma$ be such that 
  $\|Q_e\|_\infty\leq\gamma$. Then, for all $t$ it holds that
  \begin{align*}
    \tr \phaar^K \leq \tr((Q_S^R+\gamma\ii)^{2t}). 
  \end{align*}
\end{proposition}
\begin{proof}
  Let $\{r_i\}$ be the eigenvalues of $P_S^K$. By Weyl's perturbation theorem, for each $r_i$, there is some eigenvalue $q_i$ of $Q_S^R$ satisfying $q_i \in r_i \pm \gamma$. 
  In particular, $Q_S^R+\gamma\ii$ has $\tr \phaar^K$-many eigenvalues in the range $[1,1+2\gamma]$. 
  Then,
  \begin{align*}
    \tr((Q_S^R+\gamma\ii)^{2t}) 
    \geq 
    \tr \phaar^K + \sum_{\substack{i\;\mathrm{s.t.}\\r_i<1}} (q_i+\gamma)^{2t}
    \geq
    \tr \phaar^K.
  \end{align*}
\end{proof}

We will show that $\|Q_e\|_\infty\leq c_2$ in Prop.~\ref{prop:qe} from App.~\ref{app:proofs}, and so we use $\gamma=c_2$ henceforth. 
Then, if for any $t$ it holds that 
\begin{align*}
  \tr((Q_S^R+c_2\ii)^{2t})<2,
\end{align*}
$K$ is irreducible. We may expand
\begin{align}
  \tr((Q_S^R+c_2\ii)^{2t}) 
  &= 
  \sum_{k=0}^{2t} \binom{2t}{k} c_2^{2t-k}\tr ((Q_S^R)^k)\\
  &=
  \sum_{k=0}^{2t} \binom{2t}{k} c_2^{2t-k}\;
  \frac{1}{|S|^k}\sum_{\mathbf{s}\in S^k}
  |\tr \t\rho_R(\mathbf{s})|^2,\label{eq:to estimate}
\end{align}
where we used,
\begin{align*}
  \t\rho_R(s) := \t\pi_R \t\rho(s) \t\pi_R^\dagger, \qquad 
  \t\rho_R(\mbf{s}):=\t\rho_R(s_1)\dots\t\rho_R(s_k), \qquad
  s\in S, \,\mbf{s}\in S^k.
\end{align*}

Our approach is to bound the norm of all terms with $k<2t$ and estimate the one with $k=2t$. 
This is because in the regime of interest $c_2$ is small, and so terms with non-trivial powers of $c_2$ are of subleading order. 
The following proposition will be used to bound the size of subleading terms.

\begin{proposition}
  \label{prop:not too large}
  Let $R$ hold an invariance certificate with precision $\epsilon<1/2$
  and let $K$ be the minimizer in eq.~\eqref{eq:minimal invariant}. Then, for any $\mathbf{s}\in S^k$, 
  it holds that 
  \begin{align*}
    |\tr \t\rho_R(\mathbf{s})|^2
    \leq 
    \dim^2 K + d_k.
  \end{align*}
\end{proposition}

The following proposition uses the previous result to bound the size of the subleading contributions to eq.~\eqref{eq:to estimate}.

\begin{proposition}
  \label{prop:subleading}
  Let $R$, $K$ and $\epsilon$ be as in Prop.~\ref{prop:not too large}, and let $n\floating<1/2$.
  Then,
  \begin{align*}
    \left|\sum_{k=0}^{2t-1} \binom{2t}{k} c_2^{2t-k}\tr((Q_S^R)^k)\right|
    \leq
    e_t.
  \end{align*}
\end{proposition}

We therefore obtain
\begin{align*}
  \tr\phaar^K \leq e_t + \tr((Q_S^R)^{2t}) 
  = 
  e_t + \frac{1}{|S|^{2t}}\sum_{\mbf{s}\in S^{2t}}|\tr \t\rho_R(\mbf{s})|^2.
\end{align*}
All that is left to be shown is that the estimator for the second term used by Alg.~\ref{alg:noisy irrep} concentrates sharply around its mean.
For this we will use the following proposition, a simple consequence of the Chernoff bound.

\begin{proposition}
  \label{prop:chernoff}
  Let $R$, $K$ and $\epsilon$ be as in Prop.~\ref{prop:not too large}, and assume that $K$ is \emph{reducible}.
  Let $\{\mathbf{s}_i\}$ be $m$ uniformly random samples from $S^{2t}$. 
  Then, for any $\theta\in(0,1)$, it holds that
  \begin{align*}
    \mathrm{Pr}\left[
    \frac1m\sum_{i=1}^m |\tr \t\rho_R(\mathbf{s}_i)|^2 \leq (1-\theta) \tr((Q_S^R)^{2t})
    \right]
    < \exp(\frac{-\theta^2 m(2-e_t)}{2(\dim^2 K + d_{2t})}).
  \end{align*}
\end{proposition}

We may now prove the first main result of this subsection.

\begin{proof}[Proof of Thm.~\ref{thm:irrep certificate}]
  By our assumption on $m$, it holds that $\theta_m<1$.
  But then using Prop.~\ref{prop:chernoff} with $\theta=\theta_m$,
  \begin{align*}
    \mathrm{Pr}\left[
    \frac1m\sum_i |\tr \t\rho_R(\mathbf{s}_i)|^2+e_t \leq 2(1-\theta_m)
    \right]
    &\leq
    \mathrm{Pr}\left[
    \frac1m\sum_i |\tr \t\rho_R(\mathbf{s}_i)|^2+e_t \leq (1-\theta_m)\left[\tr((Q_S^R)^{2t})+e_t\right]
    \right]\\
    &\leq
    \mathrm{Pr}\left[
    \frac1m\sum_i |\tr \t\rho_R(\mathbf{s}_i)|^2 \leq (1-\theta_m) \tr((Q_S^R)^{2t})
    \right]\\
    &<
    \exp(\frac{-\theta_m^2 m (2-e_t)}{2(\dim^2 K + d_{2t})})
    <
    \pthresh.
  \end{align*}
  
\end{proof}

\section{Time Complexity}
\label{sec:time}

Here we analyse the runtime of the certification procedures proposed and discuss several ways to optimize it. 

Alg.~\ref{alg:invariance} runs in $O(n^3\log n)$ steps: the main sources of complexity are the $r=O(\log n)$ matrix products and the spectral norm appearing in line~3. 
The latter has complexity at most $O(n^3)$ through the singular value decomposition.

In practice, this last step step is significantly cheaper.
Ref.~\cite{magdon} estimates the spectral norm in time $O(n^2\log n)$.
Note that the method of~\cite{magdon} is probabilistic and so it raises the false positive rate, albeit in a controllable way.
Alternatively, the spectral norm can be bounded by the Frobenius norm in $O(n^2)$ operations.

To compute the runtime of Alg.~\ref{alg:noisy irrep} we assume that $\floating$ and $\epsilon$ are small enough that $d_{2(2+\log_2 d)}$ and $e_{2+\log_2 d}$ are non-increasing functions of $d:=\dim R$ and $n$.
Here, $d_t$ and $e_t$ are defined as in the top of Sec.~\ref{sec:noisy irrep} and we use $t=2+\log d$.
For this it is sufficient to take
\begin{align}
  \label{eq:bound epsilons}
  \epsilon < \frac{1}{48(d^2+1)(2+\log_2 d)},\qquad
  \floating  < \frac{1}{120n(d^2+1)(2+\log_2 d)}.
\end{align}
In this regime the runtime of the algorithm, as it is written in the main text, is
\begin{align}
  \label{eq:irrep scaling}
  O\left(
    n^3 d^2 \log d \left(\log \frac1\pthresh+\log\frac{1}{\conf-\pthresh}\right)
  \right). 
\end{align}
Because the false negative rate is of secondary importance for our certificate, a convenient choice is $\conf=2\pthresh$ where both terms above have the same scaling.

The main bottleneck of~\eqref{eq:irrep scaling} is the $n^3$ factor, coming from the fact that the algorithm evaluates matrix products on $\cc^{n\times n}$. 
This can be significantly reduced by either: \emph{1.}\ taking products in the group and then obtaining the image, or \emph{2.}\ restricting matrices $\t\rho_R(s)$ to the subspace $R$ first, and taking products in this smaller space. 
Letting $D$ denote the runtime of whichever of these two is faster, the runtime becomes $O(D d^2 \log d \log \pthresh^{-1})$.


\begin{acknowledgments}
  We thank Markus Heinrich and Frank Vallentin for insightful conversations.
  
  This work has been supported by the DFG (SPP1798 CoSIP), Germany’s Excellence Strategy – Cluster of Excellence Matter and Light for Quantum Computing (ML4Q) EXC2004/1, Cologne’s Key Profile Area Quantum Matter and Materials, the European Union’s Horizon 2020 research and innovation programme under the Marie Skłodowska-Curie agreement No 764759, and by the Perimeter Institute for Theoretical Physics. 
  Research at Perimeter Institute is supported in part by the Government of Canada through the Department of Innovation, Science and Economic Development Canada and by the Province of Ontario through the Ministry of Economic Development, Job Creation and Trade. 
  This publication was made possible through the support of a grant from the John Templeton Foundation. 
  The opinions expressed in this publication are those of the authors and do not necessarily reflect the views of the John Templeton Foundation.
\end{acknowledgments}


\bibliographystyle{ieeetr}
\bibliography{Refs}

\appendix


\section{Proofs}
\label{app:proofs}

\begin{proposition}
  \label{prop:qe}
  Let $Q_e$ be as in Prop.~\ref{prop:noisy over estimator}, and $c_2$ be as in the beginning of Sec.~\ref{sec:noisy irrep}. 
  Then $\|Q_e\|_\infty\leq c_2$.
\end{proposition}

\begin{proof}
  Let $\rho_K(s):= \pi_K\rho(s)\pi_K$ and $D(s):= \t\rho_R(s) - \rho_K(s)$. 
  Using subadditivity, we bound
  \begin{align*}
    \|Q_e\|_\infty 
    \leq \max_s \| 
    D(s)\otimes\bar{\rho}_K(s) + \rho_K(s)\otimes \bar{D}(s) + D(s)\otimes\bar{D}(s)
    \|_\infty
    \leq \max_s (2\|D(s)\|_\infty + \|D(s)\|_\infty^2).
  \end{align*}
  Further writing $\Delta :=\t\pi_R - \pi_K$ and $\Delta(s):=\t\rho(s)-\rho(s)$, we observe that
  \begin{align*}
    D(s) = 
    \Delta\rho(s)(\pi_K+\Delta)^\dagger + (\pi_K+\Delta)\rho(s)\Delta^\dagger 
    +
    (\pi_K+\Delta)\Delta(s)(\pi_K+\Delta)^\dagger,
  \end{align*}
  and so,
  \begin{align*}
    \| D(s)\|_\infty
    \leq
    2\|\Delta\|_\infty(1+\|\Delta\|_\infty) + \|\Delta(s)\|_\infty(1+\|\Delta\|_\infty)^2.
  \end{align*}
  We can directly bound $\| \Delta(s)\|_\infty \leq n\floating$.
  Then, becaus $R$ holds an invariance certificate with precision $\epsilon$, we deduce
  \begin{align*}
    \| \Delta \|_\infty \leq n\floating + \epsilon.
  \end{align*}
  It follows that $\|D(s)\|_\infty\leq c_1$, where $c_1$ is defined as in the top of Sec.~\ref{sec:noisy irrep}, and the claim follows.
\end{proof}

\begin{proof}[Proof of Prop.~\ref{prop:not too large}]
  As in the proof of Prop.~\ref{prop:qe}, let $D(s) := \t\rho_R(s) - \rho_K(s)$. 
  For the sake of convenience, let us introduce the following notation: $B_1(s) = D(s)$, $B_0(s)=\rho_K(s)$, and for any bit string $v\in\ff_2^k$ and $\mathbf{s}\in S^{k}$,
  \begin{align*}
    B_v(\mathbf{s}) = B_{v_1}(s_1)B_{v_2}(s_2)\cdots B_{v_k}(s_k).
  \end{align*}
  Then, using submultiplicativity, subadditivity and unitary invariance we find that
  \begin{align*}
    \left|
      \tr(\t\rho_R(\mathbf{s}))
    \right|^2
    &\leq
    \sum_{v\in\ff_2^k}
    \left| \tr B_v(\mathbf{s})  \right|^2\\
    &\leq
    \sum_{v\in\ff_2^k}\anynorm{B_v(\mathbf{s})}{F}^2\\
    &\leq
    \dim^2 K +\sum_{v\neq0}\max_s\anynorm{D(s)}{F}^{\wt(v)}\\
    &\leq
    \dim^2 K + \sum_{w=1}^k \binom{k}{w}\max_s\anynorm{D(s)}{F}^{w}\\
    &\leq 
    \dim^2 K + \left(1+\max_s\anynorm{D(s)}{F}\right)^k-1,
  \end{align*}
  where $\wt(v)$ denotes the \emph{Hamming weight} of $v$. 
  Then, because $R$ holds an invariance certificate with precision $\epsilon$, we may use an argument analogous to the proof of Prop.~\ref{prop:qe} to bound $\max_s\|D(s)\|_F$ by $c_1$ (defined in the top of Sec.~\ref{sec:noisy irrep}). 
  This finalizes the proof.
\end{proof}

\begin{proof}[Proof of Prop.~\ref{prop:subleading}.]
  We begin by observing that $d_k\leq d_{2t}$ for all $k\leq 2t$, and so Prop.~\ref{prop:not too large} implies
  \begin{align*}
    \left|\sum_{k=0}^{2t-1} \binom{2t}{k} c_2^{2t-k}\tr((Q_S^R)^k)\right|
    \leq
    [(1+c_2)^{2t}-1](\dim^2 K + d_{2t}).
  \end{align*}
  Since $\epsilon<1/2$, $\dim K = \dim R$. 
  Finally, $n\floating<1/2$ implies that $\mrm{int}(\tr\t\pi_R) = \tr \pi_R = \dim R$.
\end{proof}

\begin{proof}[Proof of Prop.~\ref{prop:chernoff}.]
  By Prop.~\ref{prop:not too large}, $|\tr \t\rho_R(\mathbf{s}_i)|^2/(\dim^2 K + d_{2t})$ is a random variable in $[0,1]$, so Chernoff's bound gives
  \begin{align*}
    \mathrm{Pr}\left[
    \frac1m\sum_i |\tr \t\rho_R(\mathbf{s}_i)|^2 \leq (1-\theta) \tr((Q_S^R)^{2t})
    \right]
    < \exp(\frac{-\theta^2 m \tr((Q_S^R)^{2t})}{2(\dim^2 K + d_{2t})})
  \end{align*}
  But by Prop.~\ref{prop:noisy over estimator} $\tr((Q_S^R+ c_2\ii)^{2t})\geq\tr\phaar^K\geq 2$, and by Prop.~\ref{prop:subleading}   $\tr((Q_S^R)^{2t}) \geq 2 - e_t$, which finishes the proof.
\end{proof}


\section{Extension to a weaker scenario}
\label{app:generalization}

Here we show how to modify our algorithms to a setting in which the user has considerably less control over the group than is assumed in the main text.
To keep the the line of argument clean, we provide only short proof sketches for the claimed statements, and include these at the end of the appendix.
In the following, the Lie algebra $\lie{g}$ of $G$ is endowed with a $G$-invariant inner product $\langle\cdot,\cdot\rangle_{\lie{g}}$ and a corresponding 2-norm $\lienorm{\cdot}{g}$.

In the current setting, the user is assumed to know $\t\rho$ evaluated on a \emph{fixed} symmetric generator set $S$. 
The set $S$ and the representation $\rho$ must also satisfy two requirements.

The first is that $S$ is not too ‘ill-conditioned':
We say that $S$ is \emph{$(\delta,k)$-dense}, if for any $g\in G$ there exists a word $s_1\cdots s_k$ of length $k$ in $S$ for which
\begin{align*}
  \lienorm{\log g^{-1}s_1\cdots s_k}{g} \leq \delta.
\end{align*}
The second requirement is that the $\rho$-images of close-by group elements are also close-by.
That is, we say that $\rho$ is $q$-bounded if it holds that
\begin{align*}
  \anynorm{\mrm{d}\rho(X)}{F}\leq q\lienorm{X}{g}, \qquad \forall X\in\lie{g},
\end{align*}
where $\mrm{d}\rho$ is the representation of $\lie{g}$ corresponding to $\rho$.
In summary, we assume that the user knows some numbers $(\delta,k,q)$ such that $S$ is $(\delta,k)$-dense and $\rho$ is $q$-bounded (we say that $(G,S,\rho)$ is $(\delta,k,q)$-well conditioned).

In the case $G$ is finite, one may take $k$ to be the Cayley diameter and $q=\delta=0$. 
When $G$ is continuous, to the best of our knowledge there are no explicit generator sets $S$ known to be $(\delta,k)$-dense.
For special unitary groups, the Solovay-Kitaev theorem provides an asymptotic result: certain generator sets are $(\delta,O(\log^4\delta^{-1}))$-dense.
In the case of $\mrm{SU}(2)$, some progress towards an explicit scaling for the Solovay-Kitaev theorem has been made in~\cite{quantitativeSolovayKitaev}.

\begin{remark}
  One can modify the algorithms presented here to use a bound on the spectral gap $\|P_S-\phaar\|_\infty$ as an input instead of $(\delta,k,q)$. 
  However, such a bound is rarely known without diagonalizing $P_S$.
  While results stating the \emph{existence} of a gap exist for a variety of compact groups, these do not quantify how large it is (e.g.~\cite{sudgap,expansionSimple,gapSimpleLie}).
  Because of this, we do not present such a modification.
\end{remark}

\subsection{Invariance certificate}

The invariance certificate in this scenario is given by Alg.~\ref{alg:invariance 2}, where we use
\begin{align}
  \label{eq:f def}
  f(x) 
  = 
  2\sqrt{2\dim R}
  \left(
    xk + 2kn\floating(n\floating+1) + 2q\delta\exp(q\delta) + 2n\floating
  \right).
\end{align}

\begin{algorithm}[H]
  \caption{Modified invariance certificate}
  \label{alg:invariance 2}
  \hspace*{\algorithmicindent} \textbf{Input:}%
  \begin{itemize}
    \item $\{\t\rho(s):\,s\in S\}\subset\cc^{n\times n}$,
    \item $\delta\in(0,1)$, $k\in\nn$, $q\in\rr_+$, \Comment{$(G,S,\rho)$ is $(\delta,k,q)$-well conditioned.}
    \item $\t\pi_R\in\cc^{n\times n}$,
    \item $\epsilon \in (0,1/2)$.
  \end{itemize}
  \hspace*{\algorithmicindent} \textbf{Output:}%
  \emph{True/False}\\
  \hspace*{\algorithmicindent}
  \begin{algorithmic}[1]
    \State Let $f$ be defined as in eq.~\eqref{eq:f def}
    \If{$f(\max_{s\in S}\anynorm{[\t\rho(s),\t\pi_R]}{F})\leq\epsilon$}
      \State \textbf{Return:} \emph{True}
    \EndIf
    \State \textbf{Return:} \emph{False}
  \end{algorithmic}
\end{algorithm}

As in the main text, the key quantity to be bounded is $\|\phaar(\pi_R) - \pi_R\|_F$.
This is acheived by the following two propositions.

\begin{proposition}
  \label{prop:small commutator}
  Let $(G,S,\rho)$ be $(\delta,k,q)$-well conditioned and assume that 
  \begin{align*}
    \anynorm{[\t\rho(s), \t\pi_R]}{F}\leq c_3, \qquad\forall\, s\in S.
  \end{align*}
  Then, for all $g\in G$ we have that
  \begin{align*}
    \anynorm{[\phaar(\pi_R)-\pi_R}{F} \leq k c_3 + 2kn\floating(n\floating + 1) + 2q\delta\exp(q\delta)
    =:
    c_4(c_3).
  \end{align*}
\end{proposition}

Putting this together with Prop.~\ref{prop:far from commutant} shows that if Alg.~\ref{alg:invariance 2} returns \emph{True}, then $R$ is approximately invariant up to precision $\epsilon$.

\subsection{Irreducibility certificate}

We now move on to the irreducibility certificate. 
For simplicity we only present the procedure in the ideal case, given by Alg.~\ref{alg:ideal irrep 2}.
The certificate is in essence the same as Alg.~\ref{alg:ideal irrep}, with the prominent difference that $S$ is not sampled at the start.
The proof of Thm.~\ref{thm:ideal irrep certificate} carries over exactly to the current case showing that this algorithm's false positive rate is at most $\pthresh$.

Alg.~\ref{alg:ideal irrep 2} furthermore includes the parameter $t$ as an input (compare line~3 of Alg.~\ref{alg:ideal irrep}).
This choice is made for the sake of performance.
Specifically, in Prop.~\ref{prop:irrep converse 2} we bound the false negative rate whenever $t$ is large enough ---this is in the same spirit as Prop~\ref{prop:irrep converse}.
Here, though, the bound on $t$ is too large to be useful in many practical settings.

Rather than using Prop.~\ref{prop:irrep converse 2} to choose $t$, we have instead tested the performance of the algorithm for different values of $t$ (see~\cite{documentationTBD}).
There it is found that, for a variety of finite group representations, taking $t\gtrsim k$ is sufficient to bring the empirical false negative rate down to zero.

\begin{algorithm}[H]
  \caption{Modified ideal irreducibility certificate}
  \label{alg:ideal irrep 2}
  \hspace*{\algorithmicindent} \textbf{Input:}%
  \begin{itemize}
    \item $\{\rho_K(s):\,s\in S\}\subset\cc^{n_K\times n_K}$,
    \item $\pthresh\in(0,1)$,
    \item $t\in\nn$.
  \end{itemize}
  \hspace*{\algorithmicindent} \textbf{Output:}%
  \emph{True/False}.\\
  \hspace*{\algorithmicindent}
  \begin{algorithmic}[1]
    \State Set $m =3\lceil n_K^2\log(1/\pthresh)\rceil+1$\;
    \vspace{1mm}
    \State Set $\theta_m= n_K\sqrt{\frac{2\log(1/\pthresh)}{m}} $\;
    \vspace{1mm}
    \State Compute $E_m = \frac1m\sum_{i=1}^m|\tr\rho_K(\mathbf{s}_i)|^2$, with $\mathbf{s}_i\in S^{2t}$ sampled uniformly\;
    \vspace{1mm}
    \If{$E_m < 2(1-\theta_m)$}\;
      \State \Return \emph{True}\;
    \EndIf\;
    \State \Return \emph{False}
  \end{algorithmic}
\end{algorithm}

We thus conclude by analysing the false negative rate of Alg.~\ref{alg:ideal irrep 2}.
This probability is intimately related to the spectral gap of $P_S^K$, ---the mixing time of random walks in $S$.
Here, we show how to obtain a bound on this spectral gap from the parameters $(\delta,k,q)$.
This result follows from Ref.~\cite[Lemma~5]{varju} up to some minor technical detail.

\begin{proposition}
  \label{prop:irrep converse 2}
  There exists a constant $c_0$ such that for any compact group $G$, generator set $S\subset G$ and \emph{irreducible} representation $\rho_K$ the following holds.
  If $(G,S,\rho_K)$ is $(\delta,k,q)$-well conditioned with $\delta\leq (c_0q)^{-c_0}$, then for any
  \begin{align*}
    t \geq \frac12 \frac{\log n-1}{\log\frac{1}{1-1/|S|k^2}},
  \end{align*}
  it holds that the probability that Alg.~\ref{alg:ideal irrep 2} returns \emph{False} upon this input is at most 
  \begin{align*}
    \exp(\frac{-m}{3\dim^2 K}  \left(\frac{2-\theta_m}{1+(n-1)(1-k^{-2}|S|^{-1})^{2t}}-1\right)^2 ).
  \end{align*}
\end{proposition}

Our approach is the following. 
Ref.~\cite[Lemma~5]{varju} gives a bound on this spectral gap as a funtion of $\delta$, $k$ and a third parameter, the \emph{maximal weight length} defined by
\begin{align*}
  \max\left\{
    \anynorm{\omega}{\lie{g}^*}^2\;\Big|\; \omega \text{ weight in }\rho_K
  \right\}.
\end{align*}
The following proposition relates this quantity to our parameter $q$, which in turn allows us to obtain a bound on the mixing time in terms of $(\delta,k,q)$.

\begin{proposition}
  \label{prop:q weight}
  Let $(K,\rho_K)$ be a unitary representation of $G$ with maximal weight length
  equal to $w$. 
  Then 
  \begin{enumerate}[label=\alph*)]
    \item $\rho_K$ is $\sqrt{w\dim K}$-bounded, 
    \item if $\rho_K$ is $q$-bounded, then $q$ must satisfy $q\geq w$.
  \end{enumerate}
\end{proposition}


\subsection{Proofs}

\begin{proof}[Proof of Prop.~\ref{prop:small commutator}]
  We directly compute that for all $s\in S$
  \begin{align*}
    \anynorm{[\rho(s), \pi_R]}{F} 
    &\leq
    c_3 + 4n\floating + 2n^2\floating^2 =: c_5.
  \end{align*}
  Similarly, for any $\mbf{s}\in S^k$, 
  \begin{align*}
    \anynorm{[\rho(\mbf{s}), \pi_R]}{F} 
    &\leq
    kc_5,
  \end{align*}
  where we used the identity $[AB,C]=A[B,C] + [A,C]B$ iteratively.
  
  Now, let $g\in G$ be arbitrary. 
  By assumption, there exists a word $g_s:=s_1 \cdots s_k$ in $S$, together with an element $g_X := \exp(X)$ for which
  \begin{align*}
    &g = g_s g_X,\\
    &\lienorm{X}{g} \leq \delta.
  \end{align*}
  Subadditivity and submultiplicativity imply that
  \begin{align*}
    \anynorm{\rho(g)-\rho(g_s)}{F} 
    &= 
    \anynorm{\exp\mrm{d}\rho(X)-\ii}{F}\\
    &\leq
    \|\mrm{d}\rho(X)\|_F\exp(\|\mrm{d}\rho(X)\|_F)\\
    &\leq
    q\delta\exp(q\delta),
  \end{align*}
  and so
  \begin{align*}
    \anynorm{[\rho(g),\t\pi]}{F} \leq 2q\delta\exp(q\delta)+ kc_5 = c_4, 
    \qquad \forall\, g\in G.
  \end{align*}
  Finally, we may use the unitarity of $\rho$ to obtain
  \begin{align*}
    \|\phaar(\pi_R)-\pi_R\|_F \leq \ee_{g\sim G}\left[\|[\rho(g),\pi_R]\|_F\right],
  \end{align*}
  which proves the claim.
\end{proof}

\begin{proof}[Proof of Prop.~\ref{prop:q weight}]
  Let $\{\omega_i\}$ be the set of weights appearing in $\rho_K$, let $\omega_0$ be a weight
  in that set with maximal length (so $\anynorm{\omega_0}{\lie{g}^*}^2=w$) and let $\lie{t}$ 
  be the Lie algebra of the maximal torus in $G$. We begin by noting that because 
  $\anynorm{\cdot}{\lie{g}}$ is invariant under the adjoint $G$-action, we know that
  \begin{align*}
    \sup_{X\in\lie{g}} 
      \frac{\anynorm{\mrm{d}\rho_K(X)}{F}^2 }{\anynorm{X}{\lie{g}}^2}
    =
    \sup_{X\in\lie{t}} 
      \frac{\anynorm{\mrm{d}\rho_K(X)}{F}^2 }{\anynorm{X}{\lie{g}}^2}.
  \end{align*}
  For any $X\in\lie{t}$,
  \begin{align}
    \label{eq:rep length torus}
    \anynorm{\mrm{d}\rho_K(X)}{F}^2 
    = 
    \sum_i |\omega_i(X)|^2 
    = 
    \sum_i |\langle\omega_i^*, X\rangle_{\lie{g}}|^2,
  \end{align}
  where $\omega_i^*$ is the dual of $\omega_i$ with respect to the invariant inner  product.
  Using Cauchy-Schwartz on eq.~\eqref{eq:rep length torus} we obtain
  \begin{align*}
    \anynorm{\mrm{d}\rho_K(X)}{F}^2 
    \leq
    \anynorm{X}{\lie{g}}^2 \sum_i \anynorm{\omega_i^*}{\lie{g}}^2
    \leq
    (w\dim K)\anynorm{X}{\lie{g}}^2,
  \end{align*}
  which proves the first statement.
  
  For the second statement, let us choose $X = \omega_0^* / \anynorm{\omega_0^*}{\lie{g}}$ in
  eq.~\eqref{eq:rep length torus}. We obtain
  \begin{align*}
    \anynorm{\mrm{d}\rho_K(X)}{F}^2 
    =
    \sum_i \frac{|\langle\omega_i^*, \omega_0^*\rangle_{\lie{g}}|^2} {\anynorm{\omega_0^*}{\lie{g}}^2}
    \geq
    \anynorm{\omega_0^*}{\lie{g}}^2
    =
    w\,.
  \end{align*}
  But $\anynorm{X}{\lie{g}}=1$ so any $q\leq w$ would be inconsistent with the equation above.
\end{proof}

\begin{proof}[Proof of Prop.~\ref{prop:irrep converse 2}]
  By Prop.~\ref{prop:q weight}, the maximal weight-length $r$ of $\rho_K$ can be at most $q$.
  Consider the random walk operator $P_S$ associated to $\rho_K$ and let $\lambda$ be the spectral norm of the restriction of $P_S$ to the traceless subspace, 
  ---by the assumption that $\rho_K$ is irreducible, we know that $\lambda<1$. 
  
  Ref.~\cite[Lemma~5]{varju} implies that there exists a universal constant $c_0>0$ such that if $\delta \leq (c_0 q)^{-c_0}$, then
  \begin{align*}
    1-\lambda \geq \frac{1}{|S|k^2}.
  \end{align*}
  Hence,
  \begin{align}
    \label{eq:rand walk well conditioned}
    \tr P_S^{2t} \leq 1 + (n-1)\left(1-\frac{1}{|S|k^2}\right)^{2t}.
  \end{align}
  Then, for any $x\leq 1$, the right-hand side is smaller than $2-x$ if and only if
  \begin{align*}
    t \geq \frac12 \frac{\log\frac{n-1}{1-x}}{\log\frac{1}{1-1/|S|k^2}} =: t_x.
  \end{align*}
  Equivalently, for any $t$ given as in the assumption of the theorem, the right-hand side of~\eqref{eq:rand walk well conditioned} is at most $2-x_t$, where
  \begin{align*}
    x_t := 1 - (n-1)(1-1/|S|k^2)^{2t}
  \end{align*}
  
  The Chernoff bound implies that for any $\alpha>0$, if $\{\mathbf{s}_i\}$ are $m$ uniform
  samples from $S^{2t}$, then
  \begin{align}
    \label{eq:chernoff ideal converse}
    \mrm{Prob}\left[
      \frac1m \sum_i |\tr\rho_K(\mathbf{s}_i)| \geq (2-x_t)(1+\alpha)
    \right]
    \leq
    \exp(-\alpha^2 m/3\dim^2 K).
  \end{align}
  Consider the choice 
  \begin{align*}
    \alpha = \frac{2-\theta_m}{2-x_t}-1,
  \end{align*}
  where $\theta_m$ is as in Line~1 of Alg.\ref{alg:ideal irrep}. 
  Then, eq.~\eqref{eq:chernoff ideal converse} becomes
  \begin{align}
    \mrm{Prob}\left[
      \frac1m \sum_i |\tr\rho_K(\mathbf{s}_i)| \geq (2-x_t)(1+\alpha)
    \right]
    \leq
    \exp(\frac{-m}{3\dim^2 K}  \left(\frac{2-\theta_m}{2-x_t}-1\right)^2 ).
  \end{align}
\end{proof}

\end{document}